\documentclass{amsart}
\usepackage{amsfonts}
\usepackage{amssymb}
\usepackage{amsmath}  
\usepackage[alphabetic]{amsrefs}
\usepackage{verbatim}
\usepackage{stmaryrd}

\input xy
\xyoption{all}

\theoremstyle{plain} 
\newtheorem{theorem}{Theorem}[section]

\newtheorem{proposition}[theorem]{Proposition}
\newtheorem{corollary}[theorem]{Corollary}

\theoremstyle{definition} 
\newtheorem{example}[theorem]{Example}
\newtheorem{examples}[theorem]{Examples}
\newtheorem{definition}[theorem]{Definition}
\newtheorem{remark}[theorem]{Remark}

\newtheorem{notation}[theorem]{Notation}

\newcommand{\Ima}{\operatorname{\rm im\,}}
\newcommand{\End}[1]{\operatorname{\rm End}_{#1}}

\newcommand{\Hom}[1]{\operatorname{{\rm Hom}}_{#1}}

\newcommand{\Ext}[2]{\operatorname{\rm Ext}^{#1}_{#2}}

\newcommand{\dimv}{\underline{\dim}\,}

\newcommand{\MOD}{\mbox{{\rm mod \!}}}

\newcommand{\Mod}{\mbox{{\rm Mod \!}}}
\newcommand{\rep}{\mbox{{\rm rep \!\!}}}

\newcommand{\Rep}{\mbox{{\rm Rep \!\!}}}
\newcommand{\proj}{\operatorname{{\rm proj }}}

\newcommand{\ind}[1]{\operatorname{\rm ind}_{#1}}

\newcommand{\Gr}[1]{\mbox{{\rm Gr}}_{#1}}

\newcommand{\demo}[1]{\textsc{Proof.} #1 \hfill $\Box$ \bigskip}

\newcommand{\cA}{\mathcal{A}}

\newcommand{\cC}{\mathcal{C}}
\newcommand{\cD}{\mathcal{D}}

\newcommand{\cT}{\mathcal{T}}
\newcommand{\cU}{\mathcal{U}}
\newcommand{\cV}{\mathcal{V}}

\newcommand{\bA}{\mathbb{A}}
\newcommand{\bC}{\mathbb{C}}
\newcommand{\bN}{\mathbb{N}}
\newcommand{\bP}{\mathbb{P}}
\newcommand{\bQ}{\mathbb{Q}}
\newcommand{\bZ}{\mathbb{Z}}

\newcommand{\be}{\mathbf{e}}
\newcommand{\bg}{\mathbf{g}}

\newcommand{\bu}{\mathbf{u}}
\newcommand{\bx}{\mathbf{x}}
\newcommand{\by}{\mathbf{y}}
\newcommand{\fm}{\mathfrak{m}}

\begin{document}

\title{Cluster characters}

\author{Pierre-Guy Plamondon}

\address{Laboratoire de Math\'ematiques d'Orsay, Univ. Paris-Sud, CNRS, Univ. Paris-Saclay, 91405 Orsay, France}

\email{pierre-guy.plamondon@math.u-psud.fr}

\thanks{The author was partially supported by the French ANR grant SC3A (ANR-15-CE40-0004-01)}


\begin{abstract}
These are lecture notes from a mini-course given at the CIMPA school in Mar del Plata, Argentina, in March 2016.  The aim of the course was to introduce cluster characters for $2$-Calabi--Yau triangulated categories and present their main properties.  The notes start with the theory of $F$-polynomials of modules over finite-dimensional algebras.  Cluster categories are then introduced, before the more general setting of $2$-Calabi--Yau triangulated categories with cluster-tilting objects is defined.  Finally, cluster characters are presented, and their use in the categorification of cluster algebras is outlined.
\end{abstract}

\maketitle

\tableofcontents

\section{Introduction}

Shortly after the introduction of cluster algebras in \cite{FZ02}, links with an impressively vast number of fields of mathematics were uncovered.  Among these is the representation theory of finite-dimensional algebras, whose links to cluster algebras became apparent in, for instance, \cite{MRZ03}, \cite{BMRRT06}, \cite{CC06}, ...  

The link between representation theory and cluster algebra has proved itself to be fruitful on both sides: on the one hand, it has allowed an understanding of cluster algebras that has led to the proof of conjectures of Fomin and Zelevinsky: see for instance \cite{FK09},  \cite{DWZ09}, \cite{GLS10}, \cite{CKLP13}, ...  .   On the other hand, it has sparked many developments in representation theory, as illustrated by the introduction of the theory of $\tau$-tilting in \cite{AIR14}, the study cluster-tilted algebras and their representations initiated in \cite{BMR07} and the study of representations of certain quivers with loops in \cite{GLS14}, among other examples.

Central in the study of this link are cluster characters.  Broadly speaking, they are maps which associate to each module over certain algebras (or object in certain triangulated categories) an element in a certain cluster algebra.  They have been introduced in \cite{CC06}, and have been studied, used and generalized for instance in \cite{CK08}, \cite{CK06}, \cite{FK09}, \cite{Palu08}, \cite{Plamondon09}, \cite{Rupel15}, ... .  

The aim of these notes is to introduce cluster characters, present some of their main properties, and show how they can be used to categorify cluster algebras.

The notes are organized as follows.  In Section \ref{sect::representations}, we introduce $F$-polynomials of modules over finite-dimensional algebras.  They can be seen as a ``homology-free'' version of cluster characters.  Their definition relies heavily on representation theory of quivers and on projective varieties called submodule Grassmannians; these are introduced first.

In Section \ref{sect::clustercat}, we introduce the cluster category of an acyclic quiver.  We first recall the notion of derived category, and we focus on examples in type $A_n$.

Section \ref{sect::2CYcat} is devoted to the introduction of an abstract setting: that of $2$-Calabi--Yau triangulated categories with cluster-tilting objects.  This setting contains that of cluster categories, and is the one used in these notes to study cluster characters.

Finally, cluster characters are introduced in Section \ref{sect::cluster-characters}, together with some of their properties leading to a categorification of cluster algebras.

The notes reflect a mini-course I gave at the CIMPA school in Mar del Plata, Argentina, in March 2016.  Each section corresponds, more or less, to a one-hour lecture.  I take this opportunity to thank the organizers of the CIMPA -- ARTA V joint meeting during which this mini-course was given.

\section{Quiver representations and submodule Grassmannians}\label{sect::representations}
In this section, we define in an elementary way the notion of quiver representation, and introduce a projective variety, the \emph{submodule Grassmannian}, whose points parametrize subrepresentations of a given representation.

\subsection{Quiver representations}
Let $k$ be a field.  We are interested in studying modules over $k$-algebras and their submodules.  A convenient setting for this is that of quiver representations.  There are many textbooks dealing with the subject, for instance \cite{Pierce88} \cite{Ringel1984} \cite{ARS} \cite{ASS} and \cite{Schiffler2014}.

\begin{definition}
 A \emph{quiver} is an oriented graph.  More precisely, a quiver $Q$ is given by a $4$-tuple $(Q_0, Q_1, s,t)$, where
 \begin{itemize}
   \item $Q_0$ is a set, whose elements are called \emph{vertices};
   \item $Q_1$ is a set, whose elements are called \emph{arrows};
   \item $s,t:Q_1\to Q_0$ are two maps, which associate to each arrow its \emph{source} or its \emph{target}, respectively.
 \end{itemize}
\end{definition}

Quivers are allowed to have multiple edges, oriented cycles and even loops.  Throughout these notes quivers will be assumed to \emph{finite}, that is to say, their sets of vertices and arrows will be finite.

\begin{example}
 We will usually number the vertices of a quiver by using natural numbers, and use letters to name the arrows.  We will represent quivers as oriented graphs.  Here is an example:
   $$
     
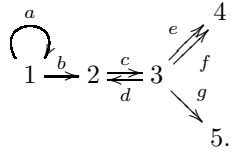
\xymatrix@!0{ &&& 4\\
                   1\ar[r]^b\ar@(ul,ur)^a & 2\ar@<.3ex>[r]^c & 3\ar@<.3ex>[l]^d\ar@<.3ex>[ur]^e\ar@<-.3ex>[ur]_f\ar[dr]^g \\
                   &&& 5. 
     }
   $$
\end{example}

A \emph{path} in a quiver is a concatenation of arrows $w=a_m\cdots a_1a_0$ such that $s(a_{i+1}) = t(a_i)$ for all $i$ from $0$ to $m-1$.  This means that we compose arrows from right to left.  We extend the maps $s$ and $t$ to the set of all paths by putting $s(w) = s(a_0)$ and $t(w) = t(a_m)$.

 Additionally, for each vertex $i$, there is a path of length $0$ starting and ending at $i$ and denoted by $e_i$.  We call it either the \emph{trivial path} or the \emph{lazy path} at $i$.  If $w$ is any path, then $w=e_{t(w)}w = we_{s(w)}$.

\begin{definition}
Let $Q$ be a quiver.  A \emph{representation of $Q$} is a tuple $V=(V_i, V_a)_{i\in Q_0, a\in Q_1}$, where
 \begin{itemize}
   \item for each vertex $i$ in $Q$, $V_i$ is a $k$-vector space, and
   \item for each arrow $a$ in $Q$, $V_a:V_{s(a)}\to V_{t(a)}$ is a $k$-linear map.
 \end{itemize}
A representation $V$ is said to be \emph{finite-dimensional} if all the vector spaces $V_i$ are finite-dimensional; in that case, the \emph{dimension vector} of $V$ is $\dimv V = (\dim V_i)_{i\in Q_0}$.  If $w=a_m\cdots a_1a_0$ is a path in $Q$, we write $V_w = V_{a_m}\circ \cdots \circ V_{a_1}\circ V_{a_0}$.
\end{definition}

 In these notes, all representations will be finite-dimensional.
 
\begin{definition}
 Let $V$ and $W$ be two representations of a quiver $Q$.  A \emph{morphism of representations from $V$ to $W$}, denoted by $f:V\to W$, is a tuple $f=(f_i)_{i\in Q_0}$, where
 \begin{itemize}
   \item for each vertex $i$ of $Q$, $f_i:V_i\to W_i$ is a $k$-linear map, and
   \item for each arrow $a$ of $Q$, we have that $W_a\circ f_{s(a)} = f_{t(a)}\circ V_a$.  In other words, the following diagram commutes:
     \begin{displaymath}
       \xymatrix{ V_{s(a)}\ar[r]^{f_{s(a)}}\ar[d]^{V_a} & W_{s(a)}\ar[d]^{W_a} \\
                  V_{t(a)}\ar[r]^{f_{t(a)}} & W_{t(a)}.
       }
     \end{displaymath}
 \end{itemize}
 Composition of morphisms is defined vertex-wise in the obvious way. 
\end{definition}

Representations of a quiver $Q$, together with their morphisms, form a category $\Rep(Q)$.  We denote by $\rep(Q)$ its full subcategory whose objects are finite-dimensional representations.  These categories are abelian; we can see this by showing that they are equivalent to module categories (see Proposition \ref{prop::Morita}).

\begin{definition}
Let $Q$ be a quiver.  The \emph{path algebra} of $Q$ is the associative $k$-algebra $kQ$ defined as follows.
  \begin{itemize}
    \item For all non-negative integers $\ell$, let $(kQ)_\ell$ be the $k$-vector with basis the set of paths of length $\ell$ in $Q$.  Then the underlying vector space of $kQ$ is $\bigoplus_{\ell=0}^{\infty} (kQ)_\ell$.   
    
    \item Multiplication is defined on paths by $$ w_2\cdot w_1 = \begin{cases} w_2w_1 & \textrm{if } s(w_2)=t(w_1) \\ 0 & \textrm{else,} \end{cases} $$
    and extended to all of $kQ$ by linearity.
  \end{itemize}
\end{definition}

We denote by $\fm$ the two-sided ideal of $kQ$ generated by the arrows of $Q$.  In other words, $\fm = \bigoplus_{\ell=1}^{\infty} (kQ)_\ell$.

If $I$ is any two-sided ideal of $kQ$, then we denote by $\Rep(Q,I)$ the full subcategory of $\Rep(Q)$ whose objects are representations $V$ ``satisfying the relations in $I$'', that is, such that for any linear combination of paths $\sum_i \lambda_iw_i$ lying in $I$, we have that $\sum_i \lambda_i V_{w_i} = 0$.

One of the main motivations for studying representations of quivers can be summarized in the following results.  First, representations of a quiver and modules over its path algebra should be viewed as being the same thing.  More precisely:

\begin{proposition}\label{prop::Morita}
  Let $Q$ be a quiver and $I$ be a two-sided ideal of $kQ$.  Then the categories $\Mod(kQ/I)$ and $\Rep(Q^{op},I^{op})$ are equivalent. (Here $\Mod A$ is the category of right(!) modules over $A$, and $Q^{op}$ is the opposite quiver, obtained by reversing the orientation of all arrows of $Q$).
  
  The same is true of $\MOD(kQ/I)$ and $\rep(Q^{op},I^{op})$, the full subcategories of finite-dimensional modules and representations, respectively.
\end{proposition}

We see $Q^{op}$ appearing in the proposition because of our choice of conventions: right modules, and composition of arrows from right to left.  The proof of the proposition is straightforward.

Secondly, over an algebraically closed field, the representation theory of \emph{any} finite-dimensional algebra is governed by a quiver with relations.  More precisely:

\begin{theorem}[Gabriel]
 Assume that the field $k$ is algebraically closed.
 For any finite-dimensional associative $k$-algebra $A$, there is a unique quiver $Q_A$ and a (non-unique) ideal $I$ of $kQ_A$ such that $A$ and $kQ_A/I$ are Morita equivalent, and $\fm^r \subset I \subset \fm^2$ for some $r\geq 2$.   
\end{theorem}

An ideal $I$ satisfying $\fm^r \subset I \subset \fm^2$ is called an \emph{admissible ideal}.

\subsection{Submodule Grassmannian}

Let $Q$ be a finite quiver, $I$ be an admissible ideal, and $V$ be a representation of $(Q,I)$.

\begin{definition}
A \emph{subrepresentation} of $V$ is a tuple $(W_i)_{i\in Q_0}$, where
 \begin{itemize}
   \item each $W_i$ is a subspace of $V_i$, and
   \item for each arrow $a$ in $Q$, we have that $V_{a}(W_{s(a)}) \subset W_{t(a)}$.
 \end{itemize}
\end{definition}
 In that case, $W= (W_i, V_a|_{W_{s(a)}})_{i\in Q_0, a\in Q_1}$ is a representation of $(Q,I)$, and the canonical inclusion into $V$ is a morphism of representations.

Grassmannians of vector spaces are projective varieties whose points parametrize subvector spaces of a given dimension.  Submodule Granssmannians of modules generalize this notion: they are projective varieties whose points parametrize submodules of a given dimension vector.

\begin{definition}\label{defi::grassmannian}
Let $\be \in \bN^{Q_0}$ be a dimension vector.  The \emph{submodule Grassmannian of $V$ with dimension vector $\be$} is the subset $\Gr{\be}(V)$ of $\prod_{i\in Q_0} \Gr{e_i}(V_i)$ of all points $(W_i)_{i\in Q_0}$ defining a subrepresentation of $V$.
\end{definition}

The submodule Grassmannian is in fact a Zariski-closed subset of $\prod_{i\in Q_0} \Gr{e_i}(V_i)$, so it is a projective variety.

\begin{examples}\label{exam::Grassmannians}
 \begin{enumerate}
  \item If the quiver $Q$ has only one vertex and no arrows, then representations of $Q$ are just vector spaces, and their submodule Grassmannians are just usual Grassmannians.
  \item\label{exam::grassKron} Let $Q =  \xymatrix{ 1\ar@<.3ex>[r] \ar@<-.3ex>[r]  & 2}$ be the Kronecker quiver.  Consider the representation
    $$
      V=\xymatrix{ k^2\ar@<.3ex>[r]^{\footnotesize\begin{pmatrix} 1 & 0 \\ 0 & 1 \end{pmatrix}} \ar@<-.3ex>[r]_{\footnotesize\begin{pmatrix} 1 & 1 \\ 0 & 1 \end{pmatrix}}  & k^2.}
    $$
    Then there are six dimension vectors for which the submodule Grassmannian of $V$ is non-empty.  The table below lists those dimension vectors and gives a variety to which the corresponding submodule Grassmannian is isomorphic.
    \[
        \begin{tabular}{|c||c|c|c|c|c|c|} 
          \hline
          $\be$   & $(0,0)$ & $(0,1)$ & $(0,2)$ & $(1,1)$ & $(1,2)$ & $(2,2)$ \\
          \hline 
          $\Gr{\be}(V)$   & point & $\bP^1$  & point & point & $\bP^1$ & point \\
          \hline
        \end{tabular}
    \]
 \end{enumerate}
\end{examples}
\subsection{$F$-polynomials of modules}

We now define the $F$-polynomial of a representation of a quiver with relations $(Q,I)$ (or, equivalently, of a module over $A=kQ/I$).  Roughly, the $F$-polynomial can be seen as a generating function for counting submodules of a given module (even though this might not make sense if the base field $k$ is infinite, since a module may have infinitely many submodules).  This theory originates from \cite{CC06}, although $F$-polynomials of modules appeared later in \cite{DWZ09}.  The general results in this section can be found in \cite[Section 2]{DG12}.

In the rest of this section, the base field $k$ is the field $\bC$ of complex numbers.

\begin{definition}\label{defi::F-polynomial}
 Let $V$ be an $A$-module. Its \emph{$F$-polynomial} is 
 $$ F_V(\by) := \sum_{\be\in\bN^{Q_0}} \chi\big( \Gr{\be}(V) \big) \by^\be,
 $$
 where
 \begin{itemize}
   \item $\by$ is the tuple of variables $(y_i \ | \ i\in Q_0)$;
   \item $\by^\be = \prod_{i\in Q_0} y_i^{e_i}$;
   \item $\Gr{\be}(V)$ is the submodule Grassmannian (see Definition \ref{defi::grassmannian}); and
   \item $\chi$ is the Euler-Poincar\'e characteristic.
 \end{itemize}
\end{definition}

We give examples of $F$-polynomials at the end of the section.  It is easy to see that the $F$-polynomial of a module $V$ only depends on the isomorphism class of $V$.

\begin{remark}\label{rema::Euler}
The most difficult part in the computation of an $F$-polynomial is determining the submodule Grassmannians $\Gr{\be}(V)$.  To compute their Euler-Poincar\'e characteristic, the following facts (true since we work over $\bC$!) are often sufficient (and indeed, suffice to prove all the formulas in these notes):
  \begin{enumerate}
    \item $\chi({\rm point})=1$;
    \item $\chi(\bA^{\! n})=1$, where $\bA^{\! n}$ is the affine space of dimension $n$;
    \item $\chi(\bP^n)=n+1$, where $\bP^n$ is the projective space of dimension $n$;
    \item $\chi(\cU \times \cV) = \chi(\cU) \cdot \chi(\cV)$;
    \item if $\cU$ is a disjoint union of two constructible subsets $C_1$ and $C_2$, then $\chi(\cU) = \chi(C_1)+\chi(C_2)$.
    \item if $f:\cU\to \cV$ is a surjective morphism of varieties (or even a surjective constructible map) such that all fibers $f^{-1}(x)$ have the same Euler characteristic, say $c$, then $\chi(\cU) = c\chi(\cV)$.
  \end{enumerate}
See \cite[Proposition 7.4.1]{GLS072}, which itself refers to \cite{MacPherson74} and \cite{Dimca04}.  On constructible maps, we refer the reader to \cite{Joyce2006}.
\end{remark}

The first property of $F$-polynomials deals with direct sums, or equivalently, with split exact sequences.

\begin{proposition}[\cite{CC06}, \cite{DG12}]\label{prop::direct-sum}
Let $V$ and $W$ be two modules over $A=kQ/I$.  Then $F_V \cdot F_W = F_{V\oplus W}$.
\end{proposition}

We outline the proof of this proposition, as it gives the flavour of the methods used to prove the various formulas that appear in these notes.  We follow \cite[Proposition 3.6]{CC06}.

{\it Proof of Proposition \ref{prop::direct-sum}}. Consider the split exact sequence 
$$ 0\to V\stackrel{\iota}{\to} V\oplus W \stackrel{\pi}{\to} W \to 0.
$$
To any submodule $B$ of $V\oplus W$ we associate the submodules $\iota^{-1}(B)$ and $\pi(B)$ of $V$ and $W$, respectively.  This defines maps
\begin{eqnarray*}
 \Phi_\be : \Gr{\be}(V\oplus W) & \longrightarrow & \coprod_{\mathbf{f} + \bg = \be} \Gr{\mathbf{f}}(V) \times \Gr{\bg}(W) \\
                 B & \longmapsto & (\iota^{-1}(B), \pi(B))
\end{eqnarray*}
which are constructible maps.  These maps are clearly surjective.  Moreover, the fiber of a point $(U_1, U_2)$ can be shown to be an affine space (it is isomorphic to $\Hom{A}(U_2, V/U_1)$, see \cite[Lemma 3.8]{CC06}).

Thus, by Remark \ref{rema::Euler}, we get
\begin{eqnarray*}
 \chi\big( \Gr{\be}(V\oplus W) \big) & = & \chi\big( \coprod_{\mathbf{f} + \bg = \be} \Gr{\mathbf{f}}(V) \times \Gr{\bg}(W) \big) \\
  & = & \sum_{\mathbf{f} + \bg = \be} \chi\big( \Gr{\mathbf{f}}(V) \big) \cdot \chi\big( \Gr{\bg}(W) \big).
\end{eqnarray*}

From there, the proof is a simple computation:
\begin{eqnarray*}
  F_V(\by)F_W(\by) & = & \Big( \sum_{\mathbf{f}\in\bN^{Q_0}} \chi\big( \Gr{\mathbf{f}}(V) \big) \by^{\mathbf{f}} \Big) \cdot \Big( \sum_{\bg\in\bN^{Q_0}} \chi\big( \Gr{\bg}(W) \big) \by^\bg \Big) \\
  & = & \sum_{\mathbf{f}, \bg} \chi\big( \Gr{\mathbf{f}}(V) \big) \chi\big( \Gr{\bg}(W) \big) \by^{\mathbf{f} + \bg} \\
  & = & \sum_\be \Big( \sum_{\mathbf{f} + \bg = \be} \chi\big( \Gr{\mathbf{f}}(V) \big) \chi\big( \Gr{\bg}(W) \big)  \Big)\by^{\be} \\
  & = & \sum_\be \chi\big( \Gr{\be}(V\oplus W) \big) \by^{\be} \\
  & = & F_{V\oplus W}(\by).
\end{eqnarray*}\hfill $\Box$ \bigskip

The second property of $F$-polynomials, and perhaps the most important one for our purposes, deals with almost-split exact sequences.  For the theory of almost-split sequences and the definition of the Auslander-Reiten translation $\tau$, we refer the reader to the notes of the courses \cite{Malicki} and \cite{Platzeck} in this volume.

\begin{theorem}[\cite{CC06}, \cite{DG12}]\label{theo::F-polynomials}
 Let $0\to \tau V\to E \to V\to 0$ be an almost-split sequence of modules over $A=kQ/I$.  Then $F_{\tau V}\cdot F_{V} = F_E + \by^{\dimv V}$.
\end{theorem}

The spirit of the proof of this theorem is similar to that of Proposition \ref{prop::direct-sum}.  The difference lies in the fact that the morphism $\Phi_{\be}$ is no longer surjective for all $\be$; the term $\by^{\dimv V}$ in the right hand side of the statement compensates, in some sense, this lack of surjectivity.

\subsection{Examples of $F$-polynomials}

\subsubsection{} Let $Q$ be the quiver with one vertex and no arrows.  Its path algebra is simply $\bC$, and representations of $Q$ are just vector spaces.  

Let $V$ be a $d$-dimensional vector space.  Then
$$ F_V(y) = \sum_{i=0}^d \binom{d}{i} y^i.
$$
This can be seen by observing that, for $d=1$, the $F$-polynomial is $1+y$, and then by applying Proposition \ref{prop::direct-sum}.  As a corollary, we get a nice proof of the known fact that the Euler-Poincar\'e characteristic of the (usual) Grassmannian $\Gr{i}(\bC^d)$ is equal to $\binom{d}{i}$.

\subsubsection{} Let $Q$ be the quiver with one vertex and one loop $\ell$, subject to the relation $\ell^2=0$.  For this quiver, there are only two indecomposable representations (up to isomorphism):
\begin{displaymath}
  \xymatrix{ V_1 = & \bC \ar@(ul,ur)^{0} & & \textrm{and} & & V_2 = & \bC^{2} \ar@(ul,ur)^{\footnotesize  \begin{pmatrix} 0 & 0 \\ 1 & 0 \end{pmatrix}}
  }
\end{displaymath}
and only one almost-split sequence:
$$ 0\to V_1\to V_2\to V_1 \to 0.
$$
The $F$-polynomials are easily seen to be $F_{V_1}(y)=1+y$ and $F_{V_2}(y)=1+y+y^2$, and one can check that they satisfy Theorem \ref{theo::F-polynomials}.

\subsubsection{} Let $Q$ and $V$ be as in Example \ref{exam::Grassmannians} (\ref{exam::grassKron}). Then $F_V(y_1, y_2)=1+2y_2+y_2^2+y_1y_2+2y_1y_2^2+y_1^2y_2^2$.

\subsubsection{}
We list a few more examples and properties of $F$-polynomials.
\begin{enumerate}
 \item If $V$ and $W$ are isomorphic, then $F_V=F_W$.  The converse is false: consider the Kronecker quiver
   $$
     \xymatrix{ 1\ar@<.3ex>[r]^a \ar@<-.3ex>[r]_b  & 2.}
   $$
  Then the representations
  $$
     \xymatrix{ V_1=\bC\ar@<.3ex>[r]^{\phantom{xx}0} \ar@<-.3ex>[r]_{\phantom{xx}1}  & \bC & \textrm{and} & V_2=\bC\ar@<.3ex>[r]^{\phantom{xx}1} \ar@<-.3ex>[r]_{\phantom{xx}0} & \bC}  
  $$
  are not isomorphic, but their $F$-polynomials are both equal to $1+y_2+y_1y_2$.
 
 \item If $F_V$ is an irreducible polynomial, then $V$ is indecomposable.  The converse is false: consider the quiver
   $$
     \xymatrix{ 1\ar@<.3ex>[r]^a   & 2\ar@<.3ex>[l]^b.}
   $$
  Then the representation
   $$
     \xymatrix{ \bC^2\ar@<.3ex>[r]^{\footnotesize\begin{pmatrix} 1 & 0 \end{pmatrix}}   & \bC\ar@<.3ex>[l]^{\footnotesize\begin{pmatrix} 0 \\ 1 \end{pmatrix}}.}
   $$
  is indecomposable, but its $F$-polynomial is $1+y_1+y_1y_2+y_1^2y_2 = (1+y_1y_2)(1+y_1)$.
 
 \item An $F$-polynomial may have negative coefficients.  An example for a quiver with two vertices and four arrows is given in \cite[Example 3.6]{DWZ09}.
\end{enumerate}

\section{Cluster categories}\label{sect::clustercat}

\subsection{Derived categories}
Derived categories were introduced by J.-L.~Verdier in \cite{Verdier77} \cite{Verdier67}.  Their general theory is discussed in numerous books and papers; let us cite \cite{Hartshorne66} \cite{Happel} \cite{Kashiwara-Schapira94} \cite{Weibel94} and \cite{Keller07}.

In this section, we only give a brief outline of the theory of derived categories, focusing on aspects that suit the purpose of these notes.  Here, $k$ is an arbitrary field.
\subsubsection{Generalities}
Let $\cA$ be an abelian category (for example, the category of modules over a finite-dimensional $k$-algebra).  In particular, every morphism in $\cA$ has a kernel and a cokernel.

A \emph{complex} of objects of $\cA$ is a sequence of morphisms
$$ C = \ldots \stackrel{d_{i-2}}{\to} C_{i-1} \stackrel{d_{i-1}}{\to} C_i \stackrel{d_{i}}{\to} C_{i+1} \stackrel{d_{i+1}}{\to} \ldots
$$
such that $d_{i+1}\circ d_i = 0$ for all integers $i$.

Let $C$ and $C'$ be two complexes.  A \emph{morphism of complexes} $f:C\to C'$ is an infinite tuple $f=(f_i)_{i \in \bZ}$ such that for all integers $i$, $f_i:C_i\to C'_i$ is a morphism, and the square
\begin{displaymath}
 \xymatrix{ C_i\ar[r]^{d_i}\ar[d]^{f_i} & C_{i+1}\ar[d]^{f_{i+1}} \\
            C'_i\ar[r]^{d'_i} & C'_{i+1}
 }
\end{displaymath}
commutes, that is, $f_{i+1}\circ d_i = d'_i\circ f_i$.

We denote by $\cC(\cA)$ the category of complexes of $\cA$.  It is an abelian category.  It admits an automorphism called the \emph{shift functor} and denoted by $[1]$, which is defined by $(C[1])_i = C_{i+1}$, and where the differential $\delta$ of $C[1]$ is defined by $\delta_i = -d_{i+1}$.

The \emph{homology} of a complex $C$ at degree $i$ is the object $H_i(C) := \ker(d_i)/\Ima(d_{i-1})$.  It is easy to see that a morphism of complexes $f:C\to C'$ induces in each degree a morphism $H_i(f):H_i(C)\to H_i(C')$.  A \emph{quasi-isomorphism} is a morphism of complexes $f$ such that all induced morphisms in homology are isomorphisms.

The derived category of $\cA$ is the category obtained when formally inverting all quasi-isomorphisms in $\cC(\cA)$.  A convenient construction of the derived category is given by first defining the \emph{homotopy category} $K(\cA)$.  This category is the quotient of $\cC(\cA)$ by the ideal of all \emph{null-homotopic morphisms}, that is, morphisms of complexes $f:C\to C'$ such that there exist morphisms $s_i:C_i\to C'_{i-1}$ is $\cA$ such that $f_i=d'_{i-1}s_i + s_{i+1}d_i$ for all $i\in \bZ$.

The \emph{derived category} $\cD(\cA)$ is then the category obtained from $K(\cA)$ by formally inverting all quasi-isomorphisms.  

If, instead of $\cC(\cA)$, one considers the categories $\cC^+(\cA)$, $\cC^-(\cA)$ and $\cC^b(\cA)$ of complexes bounded on the left, on the right and on both sides, respectively, then one defines derived categories $\cD^+(\cA)$, $\cD^-(\cA)$ and $\cD^b(\cA)$.  Of importance to us in the next section will be $\cD^b(\cA)$, called the \emph{bounded derived category}.

The advantage of defining the derived category by working in $K(\cA)$ instead of $\cC(\cA)$ is that it allows one to use a notion of ``calculus of fractions'' of morphisms, see for instance \cite[Section 2.2]{Keller07}.  Another advantage, relevant to our situation, is that if $\cA$ is the module category of a finite dimensional algebra $A$, and if we denote by $\proj A$ the full subcategory of $\cA$ whose objects are projective modules, then $\cD^-(\cA)$ is equivalent to $K^-(\proj A)$.  This latter category is often easier to work with.

\begin{proposition}
  The functor $J:\cA \to \cD^*(\cA)$ sending an object $M$ to the complex $C$ with $C_0 = M$ and $C_j=0$ if $j\neq 0$ is fully faithful.  Here, $*$ can be $+$, $-$, $b$ or an absence of symbol.
\end{proposition}

By an abuse of notation, if $M$ is an object of $\cA$, then we denote still by $M$ its image by the functor $J$.

\subsubsection{Triangulated categories}
An important property of derived categories is that they are \emph{triangulated categories}.  A triangulated category is a $k$-linear category $\cT$ together with a $k$-linear automorphism $\Sigma:\cT\to \cT$ called the \emph{suspension functor} and with a collection of sequences of morphisms of the form
$$
 X \stackrel{f}{\longrightarrow} Y \stackrel{g}{\longrightarrow} Z \stackrel{h}{\longrightarrow} \Sigma X,
$$
where $gf$ and $hg$ vanish.  The sequences belonging to the collection are called \emph{distinguished triangles}, or simply \emph{triangles}.  They are required to satisfy several axioms, which are listed below and which can be found in any of the references given at the beginning of the section.

\begin{description}
  \item[(T1)] The class of triangles is closed under isomorphism of complexes of length $4$.  For any object $X$, the sequence $X\stackrel{id_X}{\longrightarrow} X \longrightarrow 0 \longrightarrow \Sigma X$ is a triangle.  Any morphism $X \stackrel{f}{\longrightarrow} Y$ can be embedded into a triangle  $X \stackrel{f}{\longrightarrow} Y \stackrel{g}{\longrightarrow} Z \stackrel{h}{\longrightarrow} \Sigma X$.
  
  \item[(T2)] The sequence  $X \stackrel{f}{\longrightarrow} Y \stackrel{g}{\longrightarrow} Z \stackrel{h}{\longrightarrow} \Sigma X$ is a triangle if and only if $ Y \stackrel{g}{\longrightarrow} Z \stackrel{h}{\longrightarrow} \Sigma X \stackrel{-\Sigma f}{\longrightarrow} \Sigma Y$ is.

  \item[(T3)] For any commutative diagram
  $$
    \xymatrix{ X \ar[r]^f\ar[d]^u & Y\ar[r]^g\ar[d]^v & Z\ar[r]^h & \Sigma X\ar[d]^{\Sigma u} \\
               X' \ar[r]^{f'} & Y'\ar[r]^{g'} & Z'\ar[r]^{h'} & \Sigma X'
    }
  $$  
  whose rows are triangles, there exists a morphism $w:Z\longrightarrow Z'$ such that the resulting diagram also commutes (that is, $wg=g'v$ and $h'w = (\Sigma u)h$).
  
  \item[(T4)] (Octahedral axiom.) Assume that
  $$
    X \stackrel{f}{\longrightarrow} Y \stackrel{h}{\longrightarrow} Z' \stackrel{i}{\longrightarrow} \Sigma X, \quad 
    Y \stackrel{g}{\longrightarrow} Z \stackrel{j}{\longrightarrow} X' \stackrel{k}{\longrightarrow} \Sigma Y, \quad
    X \stackrel{gf}{\longrightarrow} Z \stackrel{\ell}{\longrightarrow} Y' \stackrel{m}{\longrightarrow} \Sigma X
  $$
  are triangles, and arrange them as in the following picture
  
  $\phantom{xxxx}$
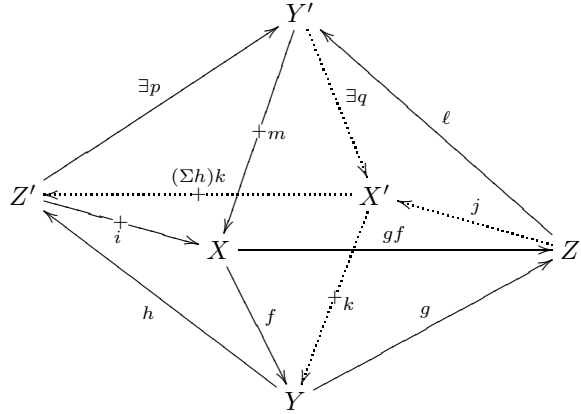
\begin{xy} 0;<1pt,0pt>:<0pt,-1pt>:: 
(105,0) *+{Y'} ="0",
(74,90) *+{X} ="1",
(207,90) *+{Z} ="2",
(0,69) *+{Z'} ="3",
(103,147) *+{Y} ="4",
(133,69) *+{X'} ="5",
"0", {\ar|+^m "1"},
"2", {\ar_\ell "0"},
"3", {\ar^{\exists p} "0"},
"0", {\ar@{.>}^{\exists q} "5"},
"1", {\ar^{gf} "2"},
"3", {\ar|+_i "1"},
"1", {\ar^f "4"},
"4", {\ar^g "2"},
"2", {\ar@{.>}_j "5"},
"4", {\ar^h "3"},
"5", {\ar@{.>}|+_{(\Sigma h)k}"3"},
"5", {\ar@{.>}^k |+"4"},
\end{xy}   

$\!\!\!\!\!\!\!$where a ``$+$'' on an arrow $A\to B$ means a morphism $A\to \Sigma B$.  Then there exist morphisms $p:Z'\to Y'$ and $q:Y'\to X'$ such that 
 $$
  Z' \stackrel{p}{\longrightarrow} Y' \stackrel{q}{\longrightarrow} X' \stackrel{(\Sigma h)k}{\longrightarrow} \Sigma Z'
 $$
 is a triangle, and we have $ph = \ell g, (\Sigma f)m = kq, i=mp$ and $j=q\ell$.
 
 In other words, the four oriented triangles in the above pictures are triangles of $\cT$, the four non-oriented triangles are commutative diagrams, and the two ``big squares'' containing the top and bottom vertices are commutative diagrams.
\end{description}

A consequence of the axioms is that for any object $X$ of $\cT$, the functor $\Hom{\cT}(X,?):\cT \to \MOD k$ sends triangles to exact sequences. 

 The derived category $\cD(\cA)$ is a triangulated category whose suspension functor is $[1]$.

\subsubsection{Hereditary case}
We now restrict to the case where $\cA = \MOD kQ$, for some finite quiver $Q$ without oriented cycles.  The path algebra $kQ$ is then \emph{hereditary}; in other words, the extension bifunctors $\Ext{i}{kQ}(?,?)$ vanish for $i\geq 2$.

In this situation, we have a good description of the objects of the bounded derived category $\cD^b(\MOD kQ)$.
\begin{proposition}[Lemma 5.2 of \cite{Happel}]
  All indecomposable objects of $\cD^b(\MOD kQ)$ are isomorphic to indecomposable \emph{stalk complexes}, that is, complexes $C$ for which there is an integer $i$ such that $C_j=0$ for $j\neq i$, and $C_i$ is an indecomposable $kQ$-module.
\end{proposition}

Thus all indecomposable objects of $\cD^b(\MOD kQ)$ have the form $M[i]$, for $M$ an indecomposable $kQ$-module and $i$ an integer.

Another important feature in this case is the existence of an automorphism of the derived category called the \emph{Auslander--Reiten translation} and denoted by $\tau$.  We refer the reader to, for instance, \cite[Section 3]{K08} for its definition in the derived category.  It is an avatar of the Auslander--Reiten translation in module categories, see \cite{Malicki} and \cite{Platzeck} in this volume, and also \cite{ARS} and \cite{ASS}.

\subsubsection{Dynkin case}
We can say even more about the structure of $\cD^b(\MOD kQ)$ if $Q$ is an orientation of a simply-laced Dynkin diagram:
$$
 \xymatrix@!0{ A_n : & 1\ar@{-}[r] & 2\ar@{-}[r] & \ldots\ar@{-}[r] & n & \\
                  & 1\ar@{-}[dr] & &&& \\
            D_n : & & 3\ar@{-}[r] & 4\ar@{-}[r] & \ldots\ar@{-}[r] & n \\
                  & 2\ar@{-}[ur] \\
                  &&& 4\ar@{-}[d] \\
            E_6 : & 1\ar@{-}[r] & 2\ar@{-}[r] & 3\ar@{-}[r] & 5\ar@{-}[r] & 6 \\
                  &&& 4\ar@{-}[d] \\
            E_7 : & 1\ar@{-}[r] & 2\ar@{-}[r] & 3\ar@{-}[r] & 5\ar@{-}[r] & 6\ar@{-}[r] & 7\\
                  &&& 4\ar@{-}[d] \\
            E_8 : & 1\ar@{-}[r] & 2\ar@{-}[r] & 3\ar@{-}[r] & 5\ar@{-}[r] & 6\ar@{-}[r] & 7\ar@{-}[r] & 8 \\
 }
$$

For any quiver $Q$, define the \emph{repetition quiver} $\bZ Q$ as follows:
\begin{itemize}
 \item vertices of $\bZ Q$ are elements $(i,n)$ of $Q_0\times \bZ$;
 \item for every arrow $a:i\to j$ in $Q$ and every integer $n$, there are arrows $(a,n):(i,n)\to (j,n)$ and $(a^*,n):(j,n)\to (i,n+1)$ in $\bZ Q$.
\end{itemize}

\begin{example}
 If $Q = 1\to 2\to 3\to 4$ is a quiver of type $A_4$, then $\bZ Q$ looks like
 $$
   
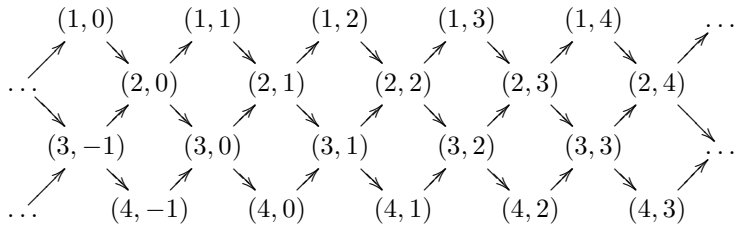
\xymatrix@!0{ &(1,0)\ar[dr] && (1,1)\ar[dr] && (1,2)\ar[dr] && (1,3)\ar[dr] && (1,4)\ar[dr] && \ldots \\
             \ldots \ar[ur]\ar[dr] && (2,0)\ar[ur]\ar[dr] && (2,1)\ar[ur]\ar[dr] && (2,2)\ar[ur]\ar[dr] && (2,3)\ar[ur]\ar[dr] && (2,4)\ar[ur]\ar[dr] & \\
              &(3,-1)\ar[ur]\ar[dr] && (3,0)\ar[ur]\ar[dr] && (3,1)\ar[ur]\ar[dr] && (3,2)\ar[ur]\ar[dr] && (3,3)\ar[ur]\ar[dr]  && \ldots &&&&\\
             \ldots \ar[ur] && (4,-1)\ar[ur] && (4,0)\ar[ur] && (4,1)\ar[ur] && (4,2)\ar[ur] && (4,3)\ar[ur]  & \\
   }
 $$
\end{example}

\begin{example}
 If $Q = \xymatrix{ 1\ar@<.3ex>[r] \ar@<-.3ex>[r]  & 2}$ is the Kronecker quiver, then $\bZ Q$ looks like
 $$
   \xymatrix@!0{& (1,0)\ar@<.3ex>[dr] \ar@<-.3ex>[dr] && (1,1)\ar@<.3ex>[dr] \ar@<-.3ex>[dr] && (1,2)\ar@<.3ex>[dr] \ar@<-.3ex>[dr] && (1,3)\ar@<.3ex>[dr] \ar@<-.3ex>[dr] && \ldots \\
             \ldots \ar@<.3ex>[ur] \ar@<-.3ex>[ur] && (2,0) \ar@<.3ex>[ur] \ar@<-.3ex>[ur] && (2,1) \ar@<.3ex>[ur] \ar@<-.3ex>[ur] && (2,2)\ar@<.3ex>[ur] \ar@<-.3ex>[ur] && (2,3)\ar@<.3ex>[ur] \ar@<-.3ex>[ur] 
   }
 $$
\end{example}

Define the \emph{mesh category} $k(\bZ Q)$ to be the category whose objects are the vertices of $\bZ Q$ and whose morphisms are $k$-linear combinations of paths in $\bZ Q$, modulo the \emph{mesh relations}: whenever we have
$$
 \xymatrix@!0{ & \bullet \ar[ddr]^{b_1} & \\
            & \bullet \ar[dr] & \\
            (i,n) \ar[uur]^{a_1}\ar[ur]\ar[dr]_{a_r} & \vdots & (i,n+1) \\
            & \bullet \ar[ur]_{b_r} & \\
 }
$$
in $\bZ Q$, where the $a_j$ are all arrows leaving $(i,n)$ and the $b_j$ are all arrows arriving in $(i,n+1)$, then $\sum_{j=1}^{r} b_ja_j = 0$.

For any category $\cC$, let $ind(\cC)$ be the full subcategory of indecomposable objects of $\cC$.

\begin{theorem}[Proposition 4.6 of \cite{Happel87}]
 If $Q$ is an orientation of a simply-laced Dynkin diagram, then $ind(\cD^b(\MOD kQ))$ is equivalent to $k(\bZ Q)$.
\end{theorem}

\subsubsection{Example: type $A_n$}\label{sect::typeA}
Many computations can be done easily in the derived category of a quiver of type $A_n$.  Let $Q=1\to 2\to \ldots \to n$.  Then $ind(\cD^b(\MOD kQ))$ is equivalent to the mesh category $k(\bZ Q)$, and so can be pictured as follows (for $n=4$):
 $${\footnotesize
   \xymatrix@!0{  &\begin{matrix} 1 \end{matrix}\ar[dr] && \begin{matrix} 2 \end{matrix}\ar[dr] && \begin{matrix} 3 \end{matrix}\ar[dr] && \begin{matrix} 4 \end{matrix}\ar[dr] && {\begin{matrix} 4\\ 3\\ 2\\ 1 \end{matrix}[1]}\ar[dr] && \ldots \\
             \ldots \ar[ur]\ar[dr] && {\begin{matrix} 2 \\ 1 \end{matrix}}\ar[ur]\ar[dr] && { \begin{matrix} 3 \\ 2 \end{matrix}}\ar[ur]\ar[dr] && {\begin{matrix} 4 \\ 3 \end{matrix}}\ar[ur]\ar[dr] &&  {\begin{matrix} 3 \\ 2 \\ 1 \end{matrix}[1]}\ar[ur]\ar[dr] && {\begin{matrix} 4 \\ 3 \\2 \end{matrix}[1]}\ar[ur]\ar[dr] & \\
              &{\begin{matrix} 4 \\ 3 \end{matrix}[-1]}\ar[ur]\ar[dr] && {\begin{matrix} 3 \\ 2 \\ 1 \end{matrix}}\ar[ur]\ar[dr] && {\begin{matrix} 4 \\ 3 \\2 \end{matrix}}\ar[ur]\ar[dr] && {\begin{matrix} 2 \\ 1 \end{matrix}[1]}\ar[ur]\ar[dr] && { \begin{matrix} 3 \\ 2 \end{matrix}[1]}\ar[ur]\ar[dr]  && \ldots &&&&\\
             \ldots \ar[ur] && \begin{matrix} 4 \end{matrix}[-1]\ar[ur] && {\begin{matrix} 4 \\ 3 \\2 \\ 1 \end{matrix}}\ar[ur] && \begin{matrix} 1 \end{matrix}[1]\ar[ur] && \begin{matrix} 2 \end{matrix}[1]\ar[ur] && \begin{matrix} 3 \end{matrix}[1]\ar[ur]  & \\
   }}
 $$
Here we denoted $kQ$-modules by their composition series (recall that right $kQ$-modules are equivalent to representations of $Q^{op}$!).  The action of the shift functor $[1]$ can be seen on the diagram; that of the Auslander-Reiten translation $\tau$ is ``translation to the left''.

Morphism spaces between two indecomposable objects can be completely determined using the mesh relations; in particular, these vector spaces have dimension at most $1$.  Some triangles can also be derived directly on the picutre: 
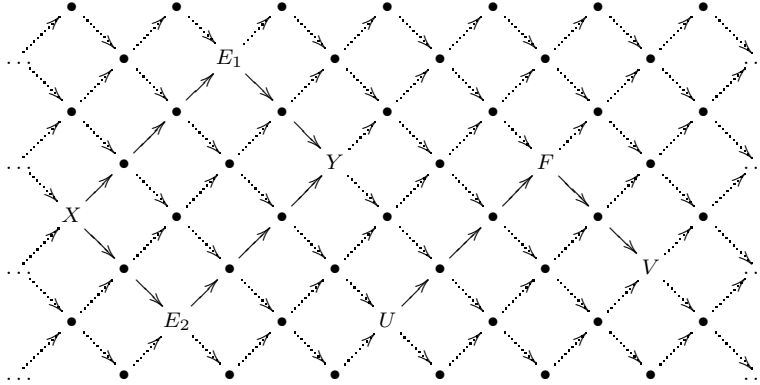
  $${\footnotesize
 \xymatrix@!0@R=0.7cm@C=0.7cm{  
             &\bullet\ar@{.>}[dr] && \bullet\ar@{.>}[dr] &&\bullet\ar@{.>}[dr] &&\bullet\ar@{.>}[dr] && \bullet\ar@{.>}[dr] &&\bullet\ar@{.>}[dr] &&\bullet\ar@{.>}[dr]  \\
             \ldots \ar@{.>}[ur]\ar@{.>}[dr] && \bullet\ar@{.>}[ur]\ar@{.>}[dr] && E_1\ar@{.>}[ur]\ar[dr] && \bullet\ar@{.>}[ur]\ar@{.>}[dr] &&  \bullet\ar@{.>}[ur]\ar@{.>}[dr] && \bullet\ar@{.>}[ur]\ar@{.>}[dr] &&\bullet\ar@{.>}[ur]\ar@{.>}[dr] &&\ldots& \\
              &\bullet\ar@{.>}[ur]\ar@{.>}[dr] &&\bullet\ar[ur]\ar@{.>}[dr] && \bullet\ar@{.>}[ur]\ar[dr] && \bullet\ar@{.>}[ur]\ar@{.>}[dr] && \bullet\ar@{.>}[ur]\ar@{.>}[dr]&&\bullet\ar@{.>}[ur]\ar@{.>}[dr] &&\bullet\ar@{.>}[ur]\ar@{.>}[dr]   &&&&\\
              \ldots \ar@{.>}[ur]\ar@{.>}[dr] && \bullet\ar[ur]\ar@{.>}[dr] && \bullet\ar@{.>}[ur]\ar@{.>}[dr] && Y\ar@{.>}[ur]\ar@{.>}[dr] &&  \bullet\ar@{.>}[ur]\ar@{.>}[dr] && F\ar@{.>}[ur]\ar[dr] &&\bullet\ar@{.>}[ur]\ar@{.>}[dr] &&\ldots& \\
              &X\ar[ur]\ar[dr] &&\bullet\ar@{.>}[ur]\ar@{.>}[dr] && \bullet\ar[ur]\ar@{.>}[dr] && \bullet\ar@{.>}[ur]\ar@{.>}[dr] && \bullet\ar[ur]\ar@{.>}[dr]&&\bullet\ar@{.>}[ur]\ar[dr] &&\bullet\ar@{.>}[ur]\ar@{.>}[dr]   &&&&\\
              \ldots \ar@{.>}[ur]\ar@{.>}[dr] && \bullet\ar@{.>}[ur]\ar[dr] && \bullet\ar[ur]\ar@{.>}[dr] && \bullet\ar@{.>}[ur]\ar@{.>}[dr] &&  \bullet\ar[ur]\ar@{.>}[dr] && \bullet\ar@{.>}[ur]\ar@{.>}[dr] &&V\ar@{.>}[ur]\ar@{.>}[dr] &&\ldots& \\
              &\bullet\ar@{.>}[ur]\ar@{.>}[dr] && E_2 \ar[ur]\ar@{.>}[dr] && \bullet\ar@{.>}[ur]\ar@{.>}[dr] && U\ar[ur]\ar@{.>}[dr] && \bullet\ar@{.>}[ur]\ar@{.>}[dr]&&\bullet\ar@{.>}[ur]\ar@{.>}[dr] &&\bullet\ar@{.>}[ur]\ar@{.>}[dr]   &&&&\\
             \ldots \ar@{.>}[ur] && \bullet\ar@{.>}[ur] &&\bullet\ar@{.>}[ur] && \bullet\ar@{.>}[ur] && \bullet\ar@{.>}[ur] && \bullet\ar@{.>}[ur] &&\bullet\ar@{.>}[ur] &&\ldots & \\
   }} 
 $$

On the left of the picture, we see a ``rectangle'' of solid arrows; it induces a triangle $X\to E_1\oplus E_2\to Y \to \Sigma X$. 

 On the right of the picture, we see a ``hook'' of solid arrows, which induces a triangle $U\to F\to V\to \Sigma U$.  The rule that ``hooks'' must obey is the following: the length of the second part of the hook (from $F$ to $V$ on the picture) is one more than the length of the downward path from the first object (here $U$) to the bottom of the picture.  Of course, hooks that are symmetric to the one pictured also yield triangles.

\subsection{Cluster categories}
Cluster categories are triangulated categories that share many of the combinatorial properties of cluster algebras.  They constitute the main setting for the definition of cluster characters (see Section \ref{sect::cluster-characters}).

\subsubsection{Orbit categories}
\begin{definition}
 Let $\cC$ be a $k$-linear category, and let $F$ be an automorphism of $\cC$.  The \emph{orbit category} $\cC/F$ is the $k$-linear category defined as follows:
 \begin{itemize}
  \item its objects are the objects of $\cC$;
  \item for any objects $X$ and $Y$, $\Hom{\cC/F}(X,Y):= \bigoplus_{n\in \bZ}\Hom{\cC}(X, F^n Y)$.
 \end{itemize}
\end{definition}
As one might expect from the name ``orbit category'', the objects $X$ and $FX$ become isomorphic in $\cC/F$.

\subsubsection{Cluster categories}
\begin{definition}[\cite{BMRRT06}]
 Let $Q$ be a quiver without oriented cycles.  The \emph{cluster category of $Q$} is the orbit category $$ \cC_Q = \cD^b(\MOD kQ)/F,$$ where $F=\tau^{-1}\circ [1]$.
\end{definition}

\begin{example}\label{exam::typeA4}
 If $Q=1\to 2\to 3\to 4$ is a quiver of type $A_4$, then using Section \ref{sect::typeA}, we get that the cluster category can be depicted as 
 $${\footnotesize
 
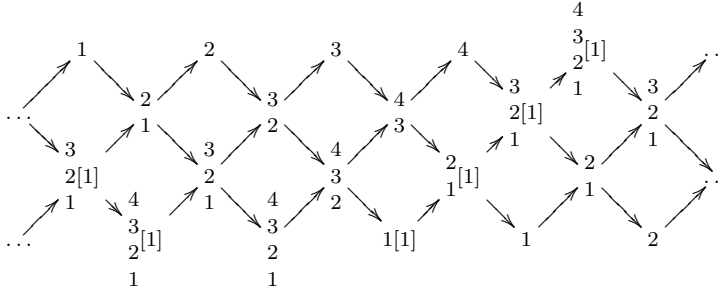
\xymatrix@!0{  &\begin{matrix} 1 \end{matrix}\ar[dr] && \begin{matrix} 2 \end{matrix}\ar[dr] && \begin{matrix} 3 \end{matrix}\ar[dr] && \begin{matrix} 4 \end{matrix}\ar[dr] && {\begin{matrix} 4\\ 3\\ 2\\ 1 \end{matrix}[1]}\ar[dr] && \ldots \\
             \ldots \ar[ur]\ar[dr] && {\begin{matrix} 2 \\ 1 \end{matrix}}\ar[ur]\ar[dr] && { \begin{matrix} 3 \\ 2 \end{matrix}}\ar[ur]\ar[dr] && {\begin{matrix} 4 \\ 3 \end{matrix}}\ar[ur]\ar[dr] &&  {\begin{matrix} 3 \\ 2 \\ 1 \end{matrix}[1]}\ar[ur]\ar[dr] && {\begin{matrix} 3 \\ 2 \\1 \end{matrix}}\ar[ur]\ar[dr] & \\
              &{\begin{matrix} 3 \\ 2 \\ 1 \end{matrix}[1]}\ar[ur]\ar[dr] && {\begin{matrix} 3 \\ 2 \\ 1 \end{matrix}}\ar[ur]\ar[dr] && {\begin{matrix} 4 \\ 3 \\2 \end{matrix}}\ar[ur]\ar[dr] && {\begin{matrix} 2 \\ 1 \end{matrix}[1]}\ar[ur]\ar[dr] && { \begin{matrix} 2 \\ 1 \end{matrix}}\ar[ur]\ar[dr]  && \ldots &&&&\\
             \ldots \ar[ur] && {\begin{matrix} 4\\ 3\\ 2\\ 1 \end{matrix}[1]}\ar[ur] && {\begin{matrix} 4 \\ 3 \\2 \\ 1 \end{matrix}}\ar[ur] && \begin{matrix} 1 \end{matrix}[1]\ar[ur] && \begin{matrix} 1 \end{matrix}\ar[ur] && \begin{matrix} 2 \end{matrix}\ar[ur]  & \\
   }} 
 $$
 Notice that the objects repeat in the diagram.  What happens is that any object $X$ becomes identified with $FX=\tau^{-1}X[1]$.  Morphism spaces and triangles can still be computed as in Section \ref{sect::typeA}.
\end{example}

Let us now list some of the most important properties of the cluster category.

\begin{theorem}[\cite{Keller05}]
 The cluster category $\cC_Q$ is a triangulated category, and the canonical functor $\cD^b(\MOD kQ)\to \cC_Q$ is a triangulated functor.
\end{theorem}

\begin{proposition}[\cite{BMR07}]
 The functor $H=\Hom{\cC_Q}(kQ, ?):\cC_Q\to \MOD kQ$ induces an equivalence of $k$-linear categories
 $$
  H:\cC_Q/(kQ[1]) \longrightarrow \MOD kQ,
 $$
 where $(kQ[1])$ is the ideal of all morphisms factoring through a direct sum of direct summands of the object $kQ[1]$.
\end{proposition}

\begin{proposition}[\cite{BMRRT06}]
 The cluster category is $2$-Calabi--Yau, in the sense of Definition \ref{defi::2CY} below.
\end{proposition}

\begin{proposition}[\cite{BMRRT06}]
 The cluster category has cluster-tilting objects, in the sense of Definition \ref{defi::CTO} below.
\end{proposition}

\section{$2$-Calabi--Yau categories}\label{sect::2CYcat}
The properties of the cluster categories listed at the end of the previous section are the ones needed for the theory of cluster characters. For this reason, we will turn to a more abstract setting where these properties are satisfied.  In this section, $k$ is an arbitrary field.
\subsection{Definition}\label{sect::2CY}
Let $\cC$ be a $k$-linear category. We will assume the following:
\begin{itemize}
 \item $\cC$ is $\Hom{}$-finite, that is, all morphism spaces in $\cC$ are finite-dimensional;
 \item $\cC$ is Krull-Schmidt, that is, every object of $\cC$ is isomorphic to a direct sum of indecomposable objects (with local endomorphism rings), and this decomposition is unique up to isomorphism and reordering of the factors;
 \item $\cC$ is triangulated, with shift functor $\Sigma$.
\end{itemize}

\begin{definition}\label{defi::2CY}
 The category $\cC$ is \emph{$2$-Calabi--Yau} if, for all objects $X$ and $Y$ of $\cC$, there is a (bifunctorial) isomorphism $$\Hom{\cC}(X,\Sigma Y) \to D\Hom{\cC}(Y,\Sigma X),$$ where $D=\Hom{k}(?,k)$ is the usual vector space duality.
\end{definition}

\begin{example}
 As seen in the previous section, the cluster category $\cC_Q$ of a quiver $Q$ without oriented cycles is a $2$-Calabi--Yau category.
\end{example}

\begin{example}
 Another family of examples is given by C.~Amiot's \emph{generalized cluster category} associated to a quiver with potential.  This is developed in \cite{Amiot08}.
\end{example}

\begin{example}
 In \cite{GLS10a} and \cite{BIRSc}, certain subcategories $\cC_w$ of the category of modules over a preprojective algebra were studied.  These categories are Frobenius categories, and their stable categories are triangulated and $2$-Calabi--Yau.
\end{example}

\subsection{Cluster-tilting objects}
Keep the notations of Section \ref{sect::2CY}.

\begin{definition}
 Let $\cC$ be a triangulated category.  An object $X$ of $\cC$ is \emph{rigid} if $\Hom{\cC}(X,\Sigma X)=0$.
\end{definition}

\begin{definition}\label{defi::CTO}
 Let $\cC$ be a $2$-Calabi--Yau category.  An object $T$ of $\cC$ is a \emph{cluster-tilting object} if
   \begin{itemize}
    \item $T$ is rigid, and
    \item for any object $X$, $\Hom{\cC}(T, \Sigma X)=0$ only if $X$ is a direct sum of direct summands of $T$.
   \end{itemize}
\end{definition}

We will usually assume that cluster-tilting objects are \emph{basic}, that is, that they can be written as a direct sum of \emph{pairwise non-isomorphic} indecomposable objects.

\begin{example}
 In a cluster category $\cC_Q$, the object $kQ$ is always a cluster-tilting object.
\end{example}

\begin{example}
An object $T$ is cluster-tilting if and only if $\Sigma T$ is.
\end{example}

\begin{example}\label{exam::CTOA4}
 In Example \ref{exam::typeA4}, the objects 
$$kQ={\footnotesize \begin{matrix} 1 \end{matrix}\oplus\begin{matrix} 2\\ 1  \end{matrix}\oplus \begin{matrix} 3\\ 2\\ 1 \end{matrix}\oplus \begin{matrix} 4\\ 3\\ 2\\ 1 \end{matrix}}, \quad T={\footnotesize \begin{matrix} 1 \end{matrix}\oplus\begin{matrix} 3  \end{matrix}\oplus \begin{matrix} 3\\ 2\\ 1 \end{matrix}\oplus \begin{matrix} 4\\ 3\\ 2\\ 1 \end{matrix}} \quad \textrm{ and }\quad T'={\footnotesize \begin{matrix} 1 \end{matrix}\oplus\begin{matrix} 3  \end{matrix}\oplus \begin{matrix} 4 \\ 3 \end{matrix}\oplus \begin{matrix} 4\\ 3\\ 2\\ 1 \end{matrix}}$$ are cluster-tilting objects.
\end{example}

We will see in Section \ref{sect::mutation} how to obtain new cluster-tilting objects from a given one.

The following property is crucial in the definition of cluster characters (Section \ref{sect::cluster-characters}): it tells us how to pass from a $2$-Calabi--Yau category to a module category.

\begin{proposition}[\cite{BMR07}\cite{KR07}]\label{prop::backtomod}
Let $T=T_1\oplus\ldots\oplus T_n$ be a basic cluster-tilting object of a $2$-Calabi--Yau category $\cC$.  We assume that the $T_i$ are indecomposable.  Then the functor 
$$
  H=\Hom{\cC}(T, \Sigma ?):\cC \longrightarrow \MOD \End{\cC}(T)
$$
induces an equivalence of $k$-linear categories
$$
 H:\cC/(T) \longrightarrow \MOD \End{\cC}(T).
$$
Moreover, 
  \begin{itemize}
   \item $H(\Sigma^{-1} T_i)$ is an indecomposable projective module for all $i\in \{1,2,\ldots, n\}$;
   \item $H(\Sigma T_i)$ is an indecomposable injective module for all $i\in \{1,2,\ldots, n\}$;
   \item for any indecomposable object $X$ other than the $T_i$, $H(\Sigma X) = \tau H(X)$, where $\tau$ is the Auslander--Reiten translation;
   \item triangles in $\cC$ are sent to long exact sequences in $\MOD \End{\cC}(T)$.
  \end{itemize}

\end{proposition}

\subsection{Index}
In a $2$-Calabi--Yau triangulated category, cluster-tilting objects act like generators of the category.  To be precise:

\begin{proposition}[\cite{KR07}]\label{prop::approximation}
 Let $\cC$ be a $2$-Calabi--Yau category with basic cluster-tilting object $T=\bigoplus_{i=1}^n T_i$.  Then for any object $X$ of $\cC$, there is a triangle
 $$
   T_1^X\to T_0^X\to X\to \Sigma T_1^X,
 $$
 where $T_0^X = \bigoplus_{i=1}^n T_i^{\oplus a_i}$ and $T_1^X=\bigoplus_{i=1}^n T_i^{\oplus b_i}$.
\end{proposition}

\begin{definition}[\cite{DK08}]
 With the notations of Proposition \ref{prop::approximation}, the \emph{index of $X$ (with respect to $T$)} is the integer vector
 $$ \ind{T}X = (a_1-b_1, \ldots, a_n-b_n).$$
\end{definition}

Note that, even though the triangle in Proposition \ref{prop::approximation} is not unique, the index is well-defined.

\begin{remark}\label{rema::indices}
 Applying $H$ to the triangle in Proposition \ref{prop::approximation}, we get an injective presentation of $HX$.  More precisely, from the triangle  
 $$
   T_1^X\to T_0^X\to X\to \Sigma T_1^X,
 $$
 we can deduce another triangle
 $$
   T_0^X\to X\to \Sigma T_1^X \to \Sigma T_0^X,
 $$
 and applying $H$ to this triangle yields the exact sequence
 $$
   0\to HX\to H(\Sigma T_1^X) \to H(\Sigma T_0^X),
 $$
 where $H(\Sigma T_0^X)$ and $H(\Sigma T_1^X)$ are injective modules by Proposition \ref{prop::backtomod}.
  This can be used to compute indices: if one can compute a minimal injective presentation of $HX$, then one can deduce the index of $X$.
\end{remark}

\begin{example}
 The index of $T_i$ is always the vector with all coordinates zero, except the $i$th one, which is $1$.  The index of $\Sigma T_i$ is the same vector, but replacing $1$ by $-1$.  These can be computed from the triangles
 $$ 0\to T_i\stackrel{id}{\to} T_i \to 0 $$ and $$T_i\to 0 \to \Sigma T_i \stackrel{id}{\to} \Sigma T_i.$$
\end{example}

\begin{example}
 Let $Q=1\to 2\to 3\to 4$, and let $\cC$ be the cluster category of $Q$, as in Example \ref{exam::typeA4}.  Take 
$$T=kQ[1] = {\footnotesize \begin{matrix} 1 \end{matrix}[1]\oplus\begin{matrix} 2\\ 1  \end{matrix}[1]\oplus \begin{matrix} 3\\ 2\\ 1 \end{matrix}[1]\oplus \begin{matrix} 4\\ 3\\ 2\\ 1 \end{matrix}[1]}. $$

Then the choice of name for the objects of $\cC$ in the figure of Example \ref{exam::typeA4} corresponds to their image by $H$ in $\MOD kQ$ (except for the summands of $T$).

We can compute the index of indecomposable objects by computing injective resolutions of modules, as pointed out in Remark \ref{rema::indices}.  The injective modules are
$$
I_1 = {\footnotesize \begin{matrix} 4\\ 3\\ 2\\ 1  \end{matrix}}, \quad I_2 = {\footnotesize \begin{matrix} 4\\ 3\\ 2 \end{matrix}}, \quad I_3 = {\footnotesize \begin{matrix} 4\\ 3  \end{matrix}} \quad \textrm{ and } \quad I_4={\footnotesize \begin{matrix} 4 \end{matrix}}.
$$

Here are some minimal injective presentations:
$$
0\to {\footnotesize \begin{matrix} 2 \end{matrix}}\to {\footnotesize \begin{matrix} 4\\ 3\\ 2 \end{matrix}}\to {\footnotesize \begin{matrix} 4\\ 3 \end{matrix}},
$$
$$
0\to {\footnotesize \begin{matrix} 3\\ 2\\ 1 \end{matrix}}\to {\footnotesize \begin{matrix} 4\\ 3\\ 2 \\ 1\end{matrix}}\to {\footnotesize \begin{matrix} 4 \end{matrix}},
$$
$$
0\to {\footnotesize \begin{matrix} 1 \end{matrix}}\to {\footnotesize \begin{matrix} 4\\ 3\\ 2\\ 1 \end{matrix}}\to {\footnotesize \begin{matrix} 4\\ 3\\ 2 \end{matrix}}.
$$
Thus $\ind{T}({\footnotesize \begin{matrix} 2 \end{matrix}})=(0,-1,1,0)$, $\ind{T}\Big({\footnotesize \begin{matrix} 3\\ 2\\ 1 \end{matrix}}\Big) = (-1,0,0,1)$ and $\ind{T}({\footnotesize \begin{matrix} 1 \end{matrix}}) = (-1,1,0,0)$.
\end{example}

Here are some properties of indices.
\begin{enumerate}
 \item\label{enum::index} For any objects $X$ and $Y$, $\ind{T}X\oplus Y = \ind{T}X + \ind{T}Y$.
 \item \cite{DK08} If $X$ and $Y$ are rigid and $\ind{T}X = \ind{T}Y$, then $X$ and $Y$ are isomorphic.
 \item \cite{Palu08} If $X\to Y\to Z\stackrel{f}{\to} \Sigma X$ is a triangle, and if $f$ lies in $(\Sigma T)$, then $\ind{T}Y = \ind{T}X+\ind{T}Z$.
 \item \cite{Palu08} For any object $X$, the vector $(\ind{T}X+\ind{T}\Sigma X)$ only depends on the dimension vector of $HX$.  
\end{enumerate}

\begin{notation}
 If $\be$ is the dimension vector of $HX$, then we put $\iota(\be):=(\ind{T}X+\ind{T}\Sigma X)$.
\end{notation}

\section{Cluster characters}\label{sect::cluster-characters}
We now come to the main aim of these notes: to define cluster characters and give some of their main properties.

In this section, $\cC$ is a $2$-Calabi--Yau category and $T=\bigoplus_{i=1}^n T_i$ is a basic cluster-tilting object of $\cC$.  The field $k$ is now assumed to be $\bC$.

\subsection{Definition}

\begin{definition}[\cite{CC06}, \cite{CK08}, \cite{Palu08}]
 The \emph{cluster character} associated to $T$ is the map $CC$ with values in $\bZ[x_1^{\pm 1}, \ldots x_n^{\pm 1}]$ defined on objects of $\cC$ by the formula
 $$ CC(X) = \bx^{\ind{T}X} \sum_{\be \in \bN^n} \chi\big( \Gr{\be}(HX) \big) \bx^{-\iota(\be)}.$$
\end{definition}

\begin{remark}
 By computing $\iota(\be)$ when $\be$ is the dimension vector of a simple module, and by using the fact that $\iota$ is additive, one can show that the above formula is equivalent to
 $$ CC(X) = \bx^{\ind{T}X} F_{HX}(\hat{y}_1, \ldots, \hat{y}_n),$$
 where 
  \begin{itemize}
    \item we define a matrix $B=(b_{ij})_{n\times n}$ by $b_{ij} = (\textrm{$\#$ arrows $i\to j$}) - (\textrm{$\#$ arrows $j\to i$})$, where arrows are taken in the Gabriel quiver of the algebra $\End{\cC}(T)$,
    \item $\hat{y}_i = \prod_{j=1}^n x_j^{b_{ji}}$, and 
    \item $F_{HX}$ is the $F$-polynomial of $HX$ as defined in Definition \ref{defi::F-polynomial}.
  \end{itemize}
\end{remark}

An immediate consequence of the above remark is the following.

\begin{proposition}[\cite{CC06}, \cite{CK08}, \cite{Palu08}]
 If $X$ and $Y$ are objects in $\cC$, then $CC(X\oplus Y) = CC(X)\cdot CC(Y)$.
\end{proposition}
\demo{This is a consequence of Proposition \ref{prop::direct-sum} and the fact that $\ind{T}X\oplus Y = \ind{T}X + \ind{T}Y$.
}

\begin{example}
For any choice of $\cC$ and $T$, we have $CC(0)=1$ and $CC(T_i)=x_i$.
\end{example}

\begin{example}
In Example \ref{exam::typeA4}, with $T=kQ[1]$, we have that $CC\Big( {\footnotesize \begin{matrix} 3 \\ 2 \\ 1 \end{matrix}} \Big) = \frac{x_1x_2 + x_1x_4 + x_3x_4 + x_2x_3x_4}{x_1x_2x_3}$.
\end{example}

\subsection{Multiplication formula}
The main theorem of the theory of cluster characters is the following \emph{multiplication formula}.

\begin{theorem}[\cite{CC06} \cite{CK08} \cite{Palu08}]\label{theo::multiplication}
 Let $X$ and $Y$ be objects of $\cC$ such that $\Hom{\cC}(X, \Sigma Y)$ is one-dimensional.  Let $\varepsilon \in \Hom{\cC}(X, \Sigma Y)$ and $\eta \in \Hom{\cC}(Y, \Sigma X)$ be non-zero (they are unique up to a scalar).  Let
 $$ Y\stackrel{i}{\to} E\stackrel{p}{\to} X\stackrel{\varepsilon}{\to} \Sigma Y$$
 and
 $$ X\stackrel{i'}{\to} E'\stackrel{p'}{\to} Y \stackrel{\eta}{\to} \Sigma X$$
 be the corresponding non-split triangles in $\cC$.  Then
 $$ CC(X)\cdot CC(Y) = CC(E) + CC(E'). $$
\end{theorem}

This result has the same spirit as Theorem \ref{theo::F-polynomials} for $F$-polynomials.  Its proof relies on the following dichotomy:

\begin{proposition}[Proposition 4.3 of \cite{Palu08}]
 Keep the notations of Theorem \ref{theo::multiplication}.  Let $U$ and $V$ be submodules of $HX$ and $HY$, respectively.  Then the two following conditions are equivalent:
 \begin{enumerate}
  \item There exists a submodule $W$ of $HE$ such that $Hp(W) = U$ and $(Hi)^{-1}(W) = V$.
  \item There does not exist any submodule $W'$ of $HE'$ such that $Hp'(W') = V$ and $(Hi')^{-1}(W')=U$.
 \end{enumerate}
\end{proposition}

This results allows us to compare Euler characteristics of the submodule Grassmannians $\Gr{*}(HE)$, $\Gr{*}(HE')$, $\Gr{*}(HX)$ and $\Gr{*}(HY)$, in a way similar to (but more involved than) what we did in the proof of Proposition \ref{prop::direct-sum}.  Together with a result concerning the indices \cite[Lemma 5.1]{Palu08}, it allows to prove Theorem \ref{theo::multiplication}.  We do not recount the proof here, but rather refer the reader to \cite{Palu08}.

\subsection{Mutation of cluster-tilting objects}\label{sect::mutation}

Assume that $R=R_1\oplus R_2\oplus \ldots \oplus R_n$ is a cluster-tilting object of $\cC$.  Assume that $\End{\cC}(R)$ is written as $\bC Q_R/I$, with $Q_R$ a finite quiver without oriented cycles of length $1$ or $2$, and $I$ an admissible ideal.

Fix $i\in \{1,\ldots, n\}$.  Consider the following triangles:
$$R_i \stackrel{\alpha}{\to} \bigoplus_{a:i\to j \textrm{ in } Q_R} R_j \to R_i^* \to \Sigma R_i$$
and
$$R^{**}_i \to \bigoplus_{b:h\to i \textrm{ in } Q_R} R_h \stackrel{\beta}{\to} R_i \to \Sigma R^{**}_i,$$
where $\alpha$ is the direct sum of all morphisms $R_i\to R_j$ corresponding to arrows $a:i\to j$ in $Q_R$, and $\beta$ is the direct sum of all morphisms $R_h\to R_i$ corresponding to arrows $b:h\to i$ in $Q_R$.

\begin{theorem}[\cite{IY08}]\label{theo::mutation}
 $\phantom{x}$
 \begin{enumerate}
  \item The objects $R_i^*$ and $R_i^{**}$ are isomorphic.
  \item The object $\mu_i(R) := R_1\oplus\ldots\oplus R_{i-1}\oplus R_i^*\oplus R_{i+1}\oplus \ldots\oplus R_n$ is a cluster-tilting object of $\cC$.
  \item The only cluster-tilting objects of $\cC$ having all $R_j$ ($j\neq i$) as direct summands are $R$ and $\mu_i(R)$.
  \item\label{one-dim} The space $\Hom{\cC}(R_i, \Sigma R_i^*)$ is one-dimensional.
 \end{enumerate}
\end{theorem}

\begin{definition}
 The object $\mu_i(R)$ of Theorem \ref{theo::mutation} is the \emph{mutation of $R$ at $i$}.  Any cluster-tilting object obtained from $R$ by a sequence of mutations is said to be \emph{reachable from $R$}.
\end{definition}

We will see in Section \ref{sect::application} why this process of mutation, coupled with the multiplication formula of Theorem \ref{theo::multiplication}, allows for a categorification of cluster algebras.

\begin{example}
In Example \ref{exam::CTOA4}, the cluster-tilting object $T$ is obtained mutating $kQ$ at the second direct summand, and $T'$ is obtained by mutating $T$ at the third direct summand.
\end{example}

An interesting result holds for cluster categories.

\begin{proposition}[Proposition 3.5 of \cite{BMRRT06}]\label{prop::graph}
 Let $Q$ be a quiver without oriented cycles, and let $T$ be a cluster-tilting object of the cluster category $\cC_Q$.  Then all cluster-tilting objects of $\cC_Q$ are reachable from $T$.
\end{proposition}

\begin{remark}
 There are $2$-Calabi--Yau categories in which cluster-tilting objects are not all reachable from each other.  An example is given in \cite[Example 4.3]{Plamondon11}.
\end{remark}

\subsection{Application: categorification of cluster algebras}\label{sect::application}
 The results of the previous sections combine neatly to provide a categorification of cluster algebras.  
 
 \begin{corollary}
  Let $T$ be a cluster-tilting object of a $2$-Calabi--Yau category $\cC$, and let $R$ be as in Theorem \ref{theo::mutation}.  Then 
  $$ CC(R_i)\cdot CC(R^*_i) = \prod_{a:i\to j \textrm{ in } Q_R} CC(R_j) + \prod_{b:h\to i \textrm{ in } Q_R} CC(R_h).$$
 \end{corollary}
 \demo{This follows directly from Theorem \ref{theo::mutation}(\ref{one-dim}) and from Theorem \ref{theo::multiplication}, and from the fact that $CC(X\oplus Y) = CC(X)\cdot CC(Y)$ for all objects $X$ and $Y$.
 }
 
 The point of this corollary is that it writes down exactly an exchange relation in a cluster algebra:
 
 \begin{definition}[\cite{FZ02}\cite{FZ03}\cite{BFZ05}\cite{FZ07}]
  Let $Q$ be a quiver with $n$ vertices without oriented cycles of length $1$ or $2$, and let $\bu=(u_1, \ldots, u_n)$ be a free generating set of the field $\bQ(x_1, \ldots, x_n)$.  Call $(Q,\bu)$ a \emph{seed}.  
  
  Then the mutation of $(Q, \bu)$ at $i$  is a new seed $(Q', \bu')$ $(u_1, \ldots, u_{i-1}, u'_{i}, u_{i+1}, \ldots, u_n)$, where
  $$ u_i\cdot u'_i = \prod_{a:i\to j \textrm{ in } Q} u_j + \prod_{b:h\to i \textrm{ in } Q} u_h.$$
 and $Q'$ is the quiver obtained from $Q$ by changing the orientation of all arrows adjacent to $i$, adding an arrow $h\to j$ for every path $h\to i \to j$, and removing cycles of length $2$.
  
 \end{definition}
 
 Now, the mutation of quivers can also be interpreted inside the cluster category:
 
 \begin{theorem}[Theorem 5.2 of \cite{BIRSm}]
 Let $R$ be as in Theorem \ref{theo::mutation}.  Assume that the endomorphism algebra of $R$ is the Jacobian algebra of a quiver with potential $(Q_R, W_R)$ (see \cite{DWZ08}).  Then the endomorphism algebra of $\mu_i(R)$ is the Jacobian algebra of the mutated quiver with potential $\mu_i(Q_R, W_R)$.  In particular, $Q_{\mu_i(R)} = \mu_i(Q_R)$.
 \end{theorem}
 
 Thus we get:
 
 \begin{corollary}[\cite{CC06}\cite{CK08}\cite{Palu08}\cite{Plamondon09}...]
  If $\cC$ is a cluster category or a generalized cluster category (see \cite{Amiot08}), then the cluster character sends reachable indecomposable objects of $\cC$ to cluster variables in the cluster algebra of $Q$, where $Q$ is the Gabriel quiver of $\End{\cC}(T)$.
 \end{corollary}

\begin{remark}
  \begin{enumerate}
  \item The multiplication formula of Theorem \ref{theo::multiplication} can be generalized to the case when the dimension of the space $\Hom{\cC}(X,\Sigma Y)$ is greater than $1$, see \cite{Palu09} and \cite{GLS072}.

  \item  Cluster characters can also be defined in the setting of stably $2$-Calabi--Yau Frobenius categories, see \cite{FK09} (and also \cite{BIRSc} and \cite{GLS10a}).

  \item  If we work over finite fields instead of $\bC$, then we can define cluster characters by counting points in submodule Grassmannians.  This leads to a categorification of quantum cluster algebras, see \cite{Rupel15}.

  \item It is possible to study cluster characters without the assumption that $\cC$ is $\Hom{}$-finite, provided $\cC$ is the generalized cluster category of a quiver with potential, see \cite{Plamondon09}.
  \end{enumerate}
\end{remark}

\begin{bibdiv}
\begin{biblist}

\bib{AIR14}{article}{
      author={Adachi, Takahide},
      author={Iyama, Osamu},
      author={Reiten, Idun},
       title={$\tau$-tilting theory},
        date={2014},
     journal={Compositio Math.},
      volume={150},
      number={3},
       pages={415\ndash 452},
}

\bib{Amiot08}{article}{
      author={Amiot, Claire},
       title={Cluster categories for algebras of global dimension 2 and quivers
  with potential},
        date={2009},
        ISSN={0373-0956},
     journal={Ann. Inst. Fourier (Grenoble)},
      volume={59},
      number={6},
       pages={2525\ndash 2590},
         url={http://aif.cedram.org/item?id=AIF_2009__59_6_2525_0},
}

\bib{ARS}{book}{
      author={Auslander, Maurice},
      author={Reiten, Idun},
      author={Smal{\o}, Sverre~O.},
       title={Representation {T}heory of {A}rtin {A}lgebras},
      series={Cambridge Studies in Advanced Mathematics},
   publisher={Cambridge University Press},
     address={Cambridge},
        date={1997},
      volume={36},
        ISBN={9780521599238},
        note={Corrected reprint of the 1995 original},
}

\bib{ASS}{book}{
      author={Assem, Ibrahim},
      author={Simson, Daniel},
      author={Skowro{\'n}ski, Andrzej},
       title={Elements of the representation theory of associative algebras.
  {V}ol. 1},
      series={London Mathematical Society Student Texts},
   publisher={Cambridge University Press},
     address={Cambridge},
        date={2006},
      volume={65},
        ISBN={978-0-521-58423-4; 978-0-521-58631-3; 0-521-58631-3},
}

\bib{BFZ05}{article}{
      author={Berenstein, Arkady},
      author={Fomin, Sergey},
      author={Zelevinsky, Andrei},
       title={Cluster algebras. {III}. {U}pper bounds and double {B}ruhat
  cells},
        date={2005},
        ISSN={0012-7094},
     journal={Duke Math. J.},
      volume={126},
      number={1},
       pages={1\ndash 52},
}

\bib{BIRSc}{article}{
      author={Buan, Aslak Bakke},
      author={Iyama, Osamu},
      author={Reiten, Idun},
      author={Scott, Jeanne},
       title={Cluster structures for 2-{C}alabi-{Y}au categories and unipotent
  groups},
        date={2009},
        ISSN={0010-437X},
     journal={Compos. Math.},
      volume={145},
      number={4},
       pages={1035\ndash 1079},
         url={http://dx.doi.org/10.1112/S0010437X09003960},
}

\bib{BIRSm}{article}{
    author={Buan, Aslak Bakke},
    author={Iyama, Osamu},
    author={Reiten, Idun},
    author={Smith, David},
     TITLE = {Mutation of cluster-tilting objects and potentials},
   JOURNAL = {Amer. J. Math.},
    VOLUME = {133},
      YEAR = {2011},
    NUMBER = {4},
     PAGES = {835\ndash 887},
}

\bib{BMR07}{article}{
      author={Buan, Aslak~Bakke},
      author={Marsh, Robert},
      author={Reiten, Idun},
       title={Cluster-tilted algebras},
        date={2007},
     journal={Trans. Amer. Math. Soc.},
      volume={359},
       pages={323\ndash 332},
}

\bib{BMRRT06}{article}{
      author={Buan, Aslak~Bakke},
      author={Marsh, Robert},
      author={Reineke, Markus},
      author={Reiten, Idun},
      author={Todorov, Gordana},
       title={Tilting theory and cluster combinatorics},
        date={2006},
        ISSN={0001-8708},
     journal={Adv. Math.},
      volume={204},
      number={2},
       pages={572\ndash 618},
         url={http://dx.doi.org/10.1016/j.aim.2005.06.003},
}

\bib{CC06}{article}{
      author={Caldero, Philippe},
      author={Chapoton, Fr{\'e}d{\'e}ric},
       title={Cluster algebras as {H}all algebras of quiver representations},
        date={2006},
        ISSN={0010-2571},
     journal={Comment. Math. Helv.},
      volume={81},
      number={3},
       pages={595\ndash 616},
         url={http://dx.doi.org/10.4171/CMH/65},
}

\bib{CKLP13}{article}{
      author={Cerulli~Irelli, Giovanni},
      author={Keller, Bernhard},
      author={Labardini-Fragoso, Daniel},
      author={Plamondon, Pierre-Guy},
       title={Linear independence of cluster monomials for skew-symmetric
  cluster algebras},
        date={2013},
     journal={Compositio Math.},
      volume={149},
      number={10},
       pages={1753\ndash 1764},
}

\bib{CK06}{article}{
      author={Caldero, Philippe},
      author={Keller, Bernhard},
       title={From triangulated categories to cluster algebras. {II}},
        date={2006},
        ISSN={0012-9593},
     journal={Ann. Sci. \'Ecole Norm. Sup. (4)},
      volume={39},
      number={6},
       pages={983\ndash 1009},
}

\bib{CK08}{article}{
      author={Caldero, Philippe},
      author={Keller, Bernhard},
       title={From triangulated categories to cluster algebras},
        date={2008},
        ISSN={0020-9910},
     journal={Invent. Math.},
      volume={172},
      number={1},
       pages={169\ndash 211},
         url={http://dx.doi.org/10.1007/s00222-008-0111-4},
}

\bib{DG12}{article}{
      author={Dominguez, Salom\'on},
      author={Geiss, Christof},
       title={A {C}aldero--{C}hapoton formula for generalized cluster
  categories},
        date={2014},
     journal={J. Algebra},
      number={399},
       pages={887\ndash 893},
}

\bib{Dimca04}{book}{
      author={Dimca, Alexandru},
       title={Sheaves in topology},
      series={Universitext},
   publisher={Springer-Verlag},
     address={Berlin},
        date={2004},
}

\bib{DK08}{article}{
      author={Dehy, Raika},
      author={Keller, Bernhard},
       title={On the combinatorics of rigid objects in 2-{C}alabi-{Y}au
  categories},
        date={2008},
        ISSN={1073-7928},
     journal={Int. Math. Res. Not. IMRN},
      number={11},
       pages={Art. ID rnn029, 17},
}

\bib{DWZ08}{article}{
   author={Derksen, Harm},
   author={Weyman, Jerzy},
   author={Zelevinsky, Andrei},
   title={Quivers with potentials and their representations. I. Mutations},
   journal={Selecta Math. (N.S.)},
   volume={14},
   date={2008},
   number={1},
   pages={59--119},
}

\bib{DWZ09}{article}{
      author={Derksen, Harm},
      author={Weyman, Jerzy},
      author={Zelevinsky, Andrei},
       title={Quivers with potentials and their representations {II}:
  applications to cluster algebras},
        date={2010},
     journal={J. Amer. Math. Soc.},
      volume={23},
      number={3},
       pages={749\ndash 790},
}

\bib{FK09}{article}{
      author={Fu, Changjian},
      author={Keller, Bernhard},
       title={On cluster algebras with coefficients and 2-{C}alabi-{Y}au
  categories},
        date={2010},
        ISSN={0002-9947},
     journal={Trans. Amer. Math. Soc.},
      volume={362},
      number={2},
       pages={859\ndash 895},
         url={http://dx.doi.org/10.1090/S0002-9947-09-04979-4},
}

\bib{FZ02}{article}{
      author={Fomin, Sergey},
      author={Zelevinsky, Andrei},
       title={Cluster algebras. {I}. {F}oundations},
        date={2002},
        ISSN={0894-0347},
     journal={J. Amer. Math. Soc.},
      volume={15},
      number={2},
       pages={497\ndash 529 (electronic)},
         url={http://dx.doi.org/10.1090/S0894-0347-01-00385-X},
}

\bib{FZ03}{article}{
      author={Fomin, Sergey},
      author={Zelevinsky, Andrei},
       title={Cluster algebras. {II}. {F}inite type classification},
        date={2003},
        ISSN={0020-9910},
     journal={Invent. Math.},
      volume={154},
      number={1},
       pages={63\ndash 121},
         url={http://dx.doi.org/10.1007/s00222-003-0302-y},
}

\bib{FZ07}{article}{
      author={Fomin, Sergey},
      author={Zelevinsky, Andrei},
       title={Cluster algebras. {IV}. {C}oefficients},
        date={2007},
        ISSN={0010-437X},
     journal={Compos. Math.},
      volume={143},
      number={1},
       pages={112\ndash 164},
         url={http://dx.doi.org/10.1112/S0010437X06002521},
}



\bib{GLS072}{article}{
      author={Geiss, Christof},
      author={Leclerc, Bernard},
      author={Schr{\"o}er, Jan},
       title={Semicanonical bases and preprojective algebras. {II}. {A}
  multiplication formula},
        date={2007},
        ISSN={0010-437X},
     journal={Compos. Math.},
      volume={143},
      number={5},
       pages={1313\ndash 1334},
}

\bib{GLS10a}{article}{
      author={Geiss, Christof},
      author={Leclerc, Bernard},
      author={Schr{\"o}er, Jan},
       title={Kac--{M}oody groups and cluster algebras},
        date={2011},
     journal={Advances in Mathematics},
      volume={228},
       pages={329\ndash 433},
}

\bib{GLS10}{article}{
      author={Geiss, Christof},
      author={Leclerc, Bernard},
      author={Schr{\"o}er, Jan},
       title={Generic bases for cluster algebras and the {C}hamber {A}nsatz},
        date={2012},
     journal={J. Amer. Math. Soc.},
      volume={25},
       pages={21\ndash 76},
}

\bib{GLS14}{article}{
      author={Geiss, Christof},
      author={Leclerc, Bernard},
      author={Schr{\"o}er, Jan},
       title={Quivers with relations for symmetrizable {C}artan matrices {I}: {F}oundations},
     journal={To appear in Inventiones Math.},
}

\bib{Happel87}{article}{
      author={Happel, Dieter},
       title={On the derived category of a finite-dimensional algebra},
        date={1987},
     journal={Comment. Math. Helv.},
      volume={62},
       pages={339\ndash 389},
}

\bib{Happel}{book}{
      author={Happel, Dieter},
       title={Triangulated categories in the representation theory of finite
  dimensional algebras},
      series={LMS Lecture Notes Series},
   publisher={Cambridge University Press},
     address={New York, New Rochelle, Melbourne, Sydney},
        date={1988},
      volume={119},
}

\bib{Hartshorne66}{book}{
      author={Hartshorne, Robin},
       title={Residues and duality},
      series={Lecture notes of a seminar on the work of A. Grothendieck, given
  at Harvard 1963/64. With an appendix by P. Deligne. Lecture Notes in
  Mathematics, No. 20},
   publisher={Springer-Verlag, Berlin-New York},
        date={1966},
}

\bib{IY08}{article}{
      author={Iyama, Osamu},
      author={Yoshino, Yuji},
       title={Mutation in triangulated categories and rigid {C}ohen-{M}acaulay
  modules},
        date={2008},
        ISSN={0020-9910},
     journal={Invent. Math.},
      volume={172},
      number={1},
       pages={117\ndash 168},
}


\bib{Joyce2006}{article}{
      author={Joyce, Dominic},
      title={Constructible functions on Artin stacks},
      date={2006},
      journal={J. London Math. Soc.},
      volume={74},
      number={3},
      pages={583\ndash 606},
}

\bib{Keller05}{article}{
      author={Keller, Bernhard},
       title={On triangulated orbit categories},
        date={2005},
     journal={Documenta Mathematica},
      volume={10},
       pages={551\ndash 581},
}

\bib{Keller07}{incollection}{
      author={Keller, Bernhard},
       title={Derived categories and tilting},
        date={2007},
   booktitle={Handbook of tilting theory},
      series={London Math. Soc. Lecture Note Ser.},
      volume={332},
   publisher={Cambridge Univ. Press, Cambridge},
       pages={49\ndash 104},
}

\bib{K08}{incollection}{
      author={Keller, Bernhard},
       title={Calabi-{Y}au triangulated categories},
        date={2008},
   booktitle={Trends in representation theory of algebras and related topics},
      series={EMS Ser. Congr. Rep.},
   publisher={Eur. Math. Soc., Z\"urich},
       pages={467\ndash 489},
         url={http://dx.doi.org/10.4171/062-1/11},
}

\bib{KR07}{article}{
      author={Keller, Bernhard},
      author={Reiten, Idun},
       title={Cluster-tilted algebras are {G}orenstein and stably
  {C}alabi-{Y}au},
        date={2007},
        ISSN={0001-8708},
     journal={Adv. Math.},
      volume={211},
      number={1},
       pages={123\ndash 151},
         url={http://dx.doi.org/10.1016/j.aim.2006.07.013},
}

\bib{Kashiwara-Schapira94}{book}{
      author={Kashiwara, Masaki},
      author={Schapira, Pierre},
       title={Sheaves on manifolds},
      series={Grundlehren der Mathematischen Wissenschaften [Fundamental
  Principles of Mathematical Sciences]},
   publisher={Springer-Verlag, Berlin},
        date={1994},
      volume={292},
        note={With a chapter in French by Christian Houzel, Corrected reprint
  of the 1990 original},
}

\bib{MacPherson74}{article}{
      author={MacPherson, Robert~D.},
       title={Chern classes for singular algebraic varieties},
        date={1974},
     journal={Ann. of Math.},
      volume={100},
      number={2},
       pages={423\ndash 432},
}

\bib{Malicki}{article}{
      author={Malicki, Piotr},
      title={Auslander-Reiten theory for finite dimensional algebras},
      journal={This volume}
}

\bib{MRZ03}{article}{
      author={Marsh, Robert},
      author={Reineke, Markus},
      author={Zelevinsky, Andrei},
       title={Generalized associahedra via quiver representations},
        date={2003},
     journal={Trans. Amer. Math. Soc.},
      volume={355},
      number={10},
       pages={4171\ndash 4186},
}

\bib{Palu08}{article}{
      author={Palu, Yann},
       title={Cluster characters for 2-{C}alabi-{Y}au triangulated categories},
        date={2008},
        ISSN={0373-0956},
     journal={Ann. Inst. Fourier (Grenoble)},
      volume={58},
      number={6},
       pages={2221\ndash 2248},
         url={http://aif.cedram.org/item?id=AIF_2008__58_6_2221_0},
}

\bib{Palu09}{article}{
      author={Palu, Yann},
       title={Cluster characters {I}{I}: {A} multiplication formula},
        date={2012},
     journal={Proc. London Math. Soc.},
      volume={104},
      number={1},
       pages={57\ndash 78},
}

\bib{Pierce88}{book}{
      author={Pierce, Richard~S.},
       title={Associative {A}lgebras},
      series={Gratuate Texts in Mathematics},
   publisher={Springer-Verlag},
     address={New York, Heidelberg, Berlin},
        date={1982},
      volume={88},
}

\bib{Plamondon09}{article}{
      author={Plamondon, Pierre-Guy},
       title={Cluster characters for cluster categories with
  infinite-dimensional morphism spaces},
        date={2011},
     journal={Advances in Mathematics},
      volume={227},
      number={1},
       pages={1 \ndash  39},
}

\bib{Plamondon11}{article}{
      author={Plamondon, Pierre-Guy},
       title={Generic bases for cluster algebras from the cluster category},
        date={2012},
     journal={Int. Math. Res. Notices},
}

\bib{Platzeck}{article}{
      author={Platzeck, Maria In\'es},
      title={Introduction to the representation theory of finite dimensional algebras},
      journal={This volume},
}

\bib{Ringel1984}{book}{
      author={Ringel, Claus~Michael},
       title={Tame algebras and integral quadratic forms},
      series={Lecture Notes in Mathematics},
   publisher={Springer-Verlag, Berlin},
        date={1984},
      volume={1099},
}

\bib{Rupel15}{article}{
      author={Rupel, Dylan},
       title={Quantum cluster characters for valued quivers},
        date={2015},
     journal={Trans. Amer. Math. Soc.},
      volume={367},
       pages={7061\ndash 7102},
}

\bib{Schiffler2014}{book}{
      author={Schiffler, Ralf},
       title={Quiver representations},
      series={CMS Books in Mathematics/Ouvrages de Math\'ematiques de la SMC},
   publisher={Springer, Cham},
        date={2014},
}

\bib{Verdier77}{inproceedings}{
      author={Verdier, Jean-Louis},
       title={Cat\'egories d\'eriv\'ees : Quelques r{\'e}sultats ({E}tat 0)},
        date={1977},
   booktitle={{C}ohomologie {E}tale: {S}{\'e}minaire de {G}{\'e}om{\'e}trie
  {A}lg{\'e}brique du {B}ois-{M}arie {SGA} 4 1/2},
   publisher={Springer Berlin Heidelberg},
     address={Berlin, Heidelberg},
       pages={262\ndash 311},
}

\bib{Verdier67}{article}{
      author={Verdier, Jean-Louis},
       title={Des cat\'egories d\'eriv\'ees des cat\'egories ab\'eliennes},
        date={1996},
     journal={Ast\'erisque},
      number={239},
       pages={xii+253 pp. (1997)},
        note={With a preface by Luc Illusie, Edited and with a note by Georges
  Maltsiniotis},
}

\bib{Weibel94}{book}{
      author={Weibel, Charles~A.},
       title={An introduction to homological algebra},
      series={Cambridge Studies in Advanced Mathematics},
   publisher={Cambridge University Press, Cambridge},
        date={1994},
      volume={38},
}

\end{biblist}
\end{bibdiv}

\end{document}